\documentclass[10pt,a4paper]{article}
\usepackage{amsmath}
\usepackage{amsfonts}
\usepackage{amssymb}
\usepackage{amsthm}
\usepackage{verbatim}
\usepackage{amscd}
\usepackage{mathrsfs}
\usepackage{longtable}
\usepackage{url}

\def\smallmat#1{\left(\begin{smallmatrix}#1\end{smallmatrix}\right)}

\usepackage[OT2,T1]{fontenc}
\DeclareSymbolFont{cyrletters}{OT2}{wncyr}{m}{n}
\DeclareMathSymbol{\Sha}{\mathalpha}{cyrletters}{"58}

\newcommand{\abs}[1]{\left\lvert#1\right\rvert}

\newcommand{\BEQ}{\begin{equation}}
\newcommand{\EEQ}{\end{equation}}

\newcommand{\C}{{\mathcal{C}}}

\newcommand{\OO}{{\mathcal{O}}}
\newcommand{\gd}{{\mathfrak d}}
\newcommand{\Res}{{\operatorname{Res}}}

\newcommand{\re}{\operatorname{Re}}

\newcommand{\Aut}{\operatorname{Aut}}

\newcommand{\Tr}{\operatorname{Tr}}
\newcommand{\Frac}{\operatorname{Frac }}

\newcommand{\im}{\operatorname{Im}}

\newcommand{\GL}{\operatorname{GL}}
\newcommand{\PGL}{\operatorname{PGL}}

\newcommand{\PSL}{\operatorname{PSL}}

\newcommand{\disc}{\operatorname{disc}}

\newcommand{\St}{\operatorname{St}}

\newcommand{\preuve}{\begin{proof}}
\newcommand{\eproof}{\end{proof}}

\newcommand{\misctitle}[1]{\par\medskip\noindent{\bf #1}.}

\newtheorem{theo}{Theorem}[section]
\newtheorem{lem}[theo]{Lemma}
\newtheorem{deff}[theo]{Definition}
\newtheorem{coro}[theo]{Corollary}
\newtheorem{prop}[theo]{Proposition}

\newtheorem*{theo*}{Theorem}
\newtheorem*{coro*}{Corollary}
\newtheorem{algo}[theo]{Algorithm}

\begin{document}
\title{An algorithm to compute relative cubic fields}
\author{Anna Morra,\\Universit\'e Rennes 1, IRMAR, \\ 263 avenue du G\'en\'eral Leclerc, CS74205, 35042 RENNES Cedex, FRANCE}
\maketitle

\begin{abstract}
\noindent Let $K$ be an imaginary quadratic number field with class number $1$.
 We describe a new, essentially linear-time algorithm,  to list all isomorphism classes of cubic extensions $L/K$ up to a bound $X$ on the norm of the relative discriminant ideal. The main tools are Taniguchi's \cite{Tan} generalization of Davenport-Heilbronn parametrisation of cubic extensions, and reduction theory for binary cubic forms over imaginary quadratic fields. Finally, we give numerical data for $K=\mathbb{Q}(i)$, and we compare our results with ray class field algorithm ones, and with asymptotic heuristics, based on a generalization of Roberts' conjecture \cite{TT}.
\end{abstract}
\section{Introduction}
Given a number field $K$, a positive integer $n$ and $X > 0$, we define
$\mathcal{F}_{K,n}(X)$ to be the set of isomorphism classes of extensions
$L/K$ such that
$$[L:K]=n,\quad \text{and}\quad\mathcal{N}_{K/\mathbb{Q}}(\gd(L/K))\leq X,$$
where $\gd(L/K)$ is the relative discriminant ideal of the
extension $L/K$. Sets of this type may be enumerated algorithmically (usually
over $\mathbb{Q}$) using the geometry of numbers, following Hunter-Martinet's
theorem~\cite{Ma2}. Asymptotically, their cardinality as $X$ tends to infinity
is the subject of folklore conjectures, predicting for instance that it should
be of the order of $X$, strikingly refined by Malle~\cite{Ma} who also fixes
the Galois group of the Galois closure of $L/K$. Small values of $n$ are of
particular interest, since computer tests become comparatively easier
and more theoretical results are available; see~\cite{Bel2} for a recent survey.

In the present paper, we will focus on the case $n = 3$.
Belabas's algorithm~\cite{Bel1} lists all representatives of
$\mathcal{F}_{\mathbb{Q},3}(X)$, in time
$O_{\varepsilon}(X^{1+\varepsilon})$, essentially linear in the size of the
output. We consider the problem of generalizing this algorithm to other base fields and we will solve it completely when $K$ is imaginary quadratic, with class number $1$. Our
main result is as follows:
\begin{theo*}
Let $K$ be an imaginary quadratic number field with class number $h_K=1$.
There exists an algorithm which lists all cubic extensions in
$\mathcal{F}_{K,3}(X)$ in time $O_\varepsilon(X^{1+\varepsilon})$, for all
$\varepsilon>0.$
\end{theo*}
\noindent For an arbitrary fixed number field $K$, Datskovsky and
Wright~\cite[Theorem I.1]{DW} proved that the cardinality of
$\mathcal{F}_{K,3}(X)$ is asymptotic to a constant (depending on $K$) times
$X$ as $X\to\infty$. It follows:

\begin{coro*}
The algorithm runs in time essentially linear in the size of the output.
\end{coro*}

The algorithm uses two main ingredients : 1) a general description of
isomorphism classes of cubic extensions $L/K$ as classes of suitable binary
quadratic forms in $K[x,y]$ modulo a $\GL_2$ action; 2) classical reduction
theory in the special case where $K$ is imaginary quadratic. Enumerating
cubic extensions then amounts to enumerating integer points in an explicit
fundamental domain, cut out by the extra condition
$\mathcal{N}_{K/\mathbb{Q}}(\gd(L/K))\leq X$.

It is interesting to compare our algorithm with the classical one, using class field theory (see Section 9.2.3 of \cite{Coh2}): the latter works in time $O_\varepsilon(X^{3/2+\varepsilon})$, unless we assume the Generalized Riemann Hypothesis to obtain $O_\varepsilon(X^{1+\varepsilon})$. So our algorithm has better \emph{unconditional} complexity. Moreover, even assuming GRH, as we did in our PARI/GP implementation, the ray class field algorithm is slower than ours (see section \ref{se:results}).

Section \ref{sec_def} is devoted to our two ingredients : Taniguchi's theorem
\cite{Tan}, which generalizes the Davenport-Heilbronn bijection used by
Belabas~\cite{Bel1}, and general facts about reduction theory for integral
binary cubic forms over imaginary quadratic fields. In Section
\ref{sec_bounds} we further assume that $K$ has class number $1$ and study
the action of $\GL_2(\mathcal{O}_K)$ on binary cubic forms and obtain a
specific fundamental domain, as well as explicit numerical bounds for the
coefficients of reduced forms.
Section \ref{sec_algo} describes the core of our algorithm and Section
\ref{sec_details} explores in detail the technical issues encountered during
the implementation of the algorithm. The final section~\ref{se:results}
presents some timings for our PARI/GP implementation, over $K=\mathbb{Q}(i)$.
\misctitle{Acknowledgments}

\noindent This work was mostly carried out during my thesis at Universit\'e Bordeaux 1, whith the support of the European Community under the Marie Curie Research Training Network GTEM (MRTN-CT-2006-035495). 

I would like to thank my advisor, Karim Belabas, for his precious help, and the Institut de Math\'ematiques de Bordeaux for the computing ressources. 

I would also like to thank John Cremona for many useful and interesting
conversations on this topic, and in particular for suggesting the contents of
section \ref{t_red}.

I would also like to thank Frank Thorne, for interesting communications about cubic fields, and comparisons of my numerical data with asymptotic results (section \ref{se:results}).

I am grateful to the anonymous referee for the useful remarks that led to this version.

Finally I would like to thank Lucia for helping me with the English corrections.

\section{Notations and preliminary results}\label{sec_def}

In this section, we recall known results, needed for our algorithm.

\subsection{Taniguchi's theorem}

\begin{deff}
Let $\mathcal{O}$ be a Dedekind domain, and let $K$ be its quotient field.

\begin{itemize}
\item Let $\mathcal{C}(\mathcal{O})$ be the set of ``cubic
    algebras'' that is, isomorphism classes of $\mathcal{O}$-algebras
  that are projective of rank $3$ as $\mathcal{O}$-modules.
\item For every fractional ideal $\mathfrak{a}$ of $\mathcal{O}$ we
  define
$$\mathcal{C}(\mathcal{O},\mathfrak{a})=\{R\in
  \mathcal{C}(\mathcal{O})\mid \St(R)=\text {the ideal class of $\mathfrak{a}$}\},$$
where $\St(R)\in Cl(\mathcal{O})$ is the Steinitz class of $R$, thus
$R$ is of the form $\omega_1\mathcal{O}\oplus
\omega_2\mathcal{O}\oplus \omega_3\mathfrak{a}$, for appropriate
$\omega_1,\omega_2,\omega_3\in \Frac(R):=R\otimes_\mathcal{O} K$. 
We define the discriminant ideal $\mathfrak{d} (R) =
\disc(\omega_1,\omega_2,\omega_3)\mathfrak{a}^2$, where as usual
$\disc(\omega_1,\omega_2,\omega_3) = \det \Tr_{\Frac(R)/K} (\omega_i\omega_j)$.

\item Let further
  $$G_{\mathfrak{a}}=\left\{\begin{pmatrix}\alpha\in\mathcal{O} & \beta\in\mathfrak{a}^{-1}\\ \gamma\in\mathfrak{a} &
  \delta\in\mathcal{O}\end{pmatrix} \Bigg\vert\quad
\alpha\delta-\beta\gamma \in\mathcal{O}^{\times}\right\},$$
$$V_{\mathfrak{a}}=\{F=(a,b,c,d) \mid
a\in\mathfrak{a},b\in\mathcal{O},c\in\mathfrak{a}^{-1},d\in\mathfrak{a}^{-2}\}.$$
If $F\in V_{\mathfrak{a}}$, its discriminant $\disc(F) = b^2c^2 - 27a^2d^2
+ 18abcd - 4ac^3 - 4b^3d$ belongs to $\mathfrak{a}^{-2}$.
\item  We consider elements of $V_{\mathfrak{a}}$ as binary cubic forms,
under the identification
$(a,b,c,d)=ax^3+bx^2y+cxy^2+dy^3$ and
we define a left-action of $G_\mathfrak{a}$ on $V_\mathfrak{a}$ by 
$$M\cdot F=(\det M)^{-1}F(\alpha x+\beta y,\gamma x+\delta y),$$
where $M=\begin{pmatrix}\alpha & \beta \\ \gamma & \delta \end{pmatrix}
\in G_\mathfrak{a}$.
\end{itemize}
\end{deff}
\noindent The following theorem generalizes the Davenport-Heilbronn~\cite{DH} theory,
corresponding to the special case $\mathcal{O} = \mathbb{Z}$, to cubic
algebras over an arbitrary Dedekind domain $\mathcal{O}$:
\begin{theo}[Taniguchi~\cite{Tan}]
There exists a canonical bijection between
$\mathcal{C}(\mathcal{O},\mathfrak{a})$ and
$V_{\mathfrak{a}}/G_{\mathfrak{a}}$ such that the following diagram is
commutative:
\[
\begin{CD}
 V_\mathfrak{a} / G_\mathfrak{a} @>>>
  \mathcal{C}(\mathcal{O},\mathfrak{a})\\
 @V{\disc}VV @VV\mathfrak{d}V \\
    \mathfrak{a}^{-2} / (\mathcal{O}^{\times})^2 @>\times\mathfrak{a}^2>>
    \left\{\text{integral ideals of $\mathcal{O}$}\right\}
\end{CD},
\]
where $\mathfrak{d}$ is the relative discriminant ideal map.
\end{theo}

\misctitle{Remarks}
\begin{itemize}
\item A computation proves that the vertical ``$\text{disc}$'' is
well defined. The other vertical map $\mathfrak{d}$ is well-defined since an
$\mathcal{O}$-algebra isomorphism preserves the discriminant. 
\item We slightly changed the notation from Taniguchi's paper, to keep consistent with the notation of the following sections (Taniguchi's action $M*F$ is given by $(M^t)\cdot F$).
\end{itemize}
\begin{coro}\label{hk1}
Let $K$ be a number field with class number $h_K=1$. Let
$\mathcal{O}=\mathcal{O}_K$ be its ring of integers. Then Taniguchi's
bijection simplifies to a bijection between binary cubic forms with
coefficients in $\mathcal{O}$ modulo $\GL_2(\mathcal{O})$ and cubic
$\mathcal{O}$-algebras.
\end{coro}

To enumerate relative cubic extensions $L/K$, we shall select only the cubic
$\mathcal{O}$-algebras $R$ which are both domains and integrally closed:
those algebras are exactly the classes of the $\mathcal{O}_L$. The algebra
$R$ is a domain if and only if $F$ is irreducible over $K$. Being integrally
closed is a local property; it is equivalent to $\mathfrak{p}$-maximality at
all prime ideals $\mathfrak{p}\subset \mathcal{O}_K$ such that
$\mathfrak{p}^2\mid \mathfrak{d}(R)$ and this can be tested using
Dedekind's criterion~\cite[Theorem 2.4.8]{Coh2}. As was done in \cite{Bel1},
it is possible to use sieve methods to control the complexity of this step by
avoiding costly discriminant factorizations.

\subsection{Fundamental domains in hyperbolic $3$-space}
In this section, we describe fundamental domains for the action of Bianchi
groups on hyperbolic $3$-space, which underlie the reduction of binary
Hermitian and cubic forms (to be dealt with in the next two sections).

\begin{deff}
Let $\mathbb{H} = \mathbb{R} + \mathbb{R} i + \mathbb{R} j + \mathbb{R} k$ be
the algebra of quaternions, let $\mathbb{C} = \mathbb{R} + \mathbb{R}i$ be
the subfield of complex numbers, and let
\begin{eqnarray*}\mathcal{H}_3&=&\{z+tj\mid
z\in\mathbb{C},t\in\mathbb{R}^*_+\}\\&=&\{h=z+tj\mid h\in\mathbb{H},
\text{ such that the $k$-component is $0$},t>0 \},\end{eqnarray*}
denote hyperbolic $3$-space.
We define the action of $\PGL_2(\mathbb{C})$ on $\mathcal{H}_3$  by
$M\cdot(z+tj)=(z'+t'j)$, with
\begin{equation}\label{hyp}\left\{\begin{array}{l}
    z'=\dfrac{\rho^2 A\overline{C}+z
    A\overline{D}+\overline{z}B\overline{C}+B\overline{D}}{\rho^2\vert
    C\vert^2+z C\overline{D}+\overline{z}
    \overline{C}D+\vert
    D\vert^2}\\  t'=\dfrac{\abs{\det(M)} t}{\rho^2\vert C\vert^2+z C\overline{D}+\overline{z}
    \overline{C}D+\vert D\vert^2},\end{array}\right.\end{equation}
where $M=\begin{pmatrix}A &
    B \\ C & D \end{pmatrix}\in\PGL_2(\mathbb{C})$ and $\rho^2=\vert
z\vert^2+t^2$.
\end{deff} 

\misctitle{Remark}
With the quaternion notations (and operations), this
translates to the neater formula
$$M\cdot h=(Ah+B)(Ch+D)^{-1}.$$

\begin{deff}\label{fund_dom}
Let $K=\mathbb{Q}(\sqrt{d_K})$ be an imaginary quadratic field of 
discriminant $d_K < 0$ and class number $1$. We define
\begin{eqnarray*}
{F}_{\mathbb{Q}(i)}&=&\left\{z\in\mathbb{C}\left\lvert\quad  0\leq
 \re(z) \leq\frac{1}{2}, 0\leq \im(z)\leq \frac{1}{2} \right.\right\},\\
{F}_{\mathbb{Q}(\sqrt{-2})}&=&\left\{z\in\mathbb{C}\left\lvert\quad -1/2\leq
 \re(z) \leq\frac{1}{2}, 0\leq \im(z)\leq \frac{\sqrt{2}}{4}\right.\right\},\\
F_{\mathbb{Q}(\sqrt{-3})}&=&\left\{z\in\mathbb{C}\left\lvert
0\leq \re(z)\leq \frac{1}{2},-\frac{\sqrt{3}}{3}\re(z)\leq \im(z)\leq\frac{\sqrt{3}}{3}\re(z)\right.\right\},\\
F_K &=& \left\{z\in\mathbb{C}\left\lvert \abs{\re(z)}\leq 1/2,0\leq
\im(z)\leq \sqrt{\vert d_K\vert }/4 \right.\right\},\textrm{ when }d_K\neq -2,-3,-4.\end{eqnarray*}
Moreover we set
$$\mathcal{B}_K =\left\{ z+tj\in\mathcal{H}_3\left\lvert
z\in F_K\textrm{ and }\abs{z}^2+t^2\geq 1\right.\right\}.$$
Let
\begin{eqnarray*}\mathcal{F}_{\mathbb{Q}(i)}&=&\left\{z+tj\in\mathcal{B}_K\left\lvert\quad \re(z)\leq \im(z) \textrm{ if }z+tj \in \partial\mathcal{B}_K\right.\right\}\\
\mathcal{F}_{\mathbb{Q}(\sqrt{-2})}&=&\left\{z+tj\in\mathcal{B}_K\left\lvert\quad \re(z)\geq 0\textrm{ if }z+tj\in\partial\mathcal{B}_K \right.\right\}\\
\mathcal{F}_{\mathbb{Q}(\sqrt{-3})}&=&\left\{z+tj\in \mathcal{B}_K\left\lvert\quad \im(z)\geq 0 \textrm{ if }z+tj\in \partial\mathcal{B}_K\right.\right\},\end{eqnarray*}
where $\partial\mathcal{B}_K$ denotes the boundary of $\mathcal{B}_K$.

Finally for $K$ such that $d_K\neq -2,-3,-4$, we define
\begin{eqnarray*}\mathcal{F}_K&=&\left\{z+tj\in\mathcal{B}_K\left\lvert \begin{array}{l}\re(z)\leq P/4 \textrm{ if }\left(\im(z)=\sqrt{\abs{d_K}}/4 \textrm{ and }\abs{\re(z)}\leq \abs{z}^2+t^2-3/4\right)  \\ \re(z)\geq 0 \textrm{ if }\left(\abs{z}^2+t^2=1 \textrm{ or }\abs{\re(z)}=1/2 \textrm{ or }\im(z)=0\right)\end{array}\right.\right\}.
\end{eqnarray*}
\end{deff}

\begin{theo}\label{boundt}
  Let $K$ be an imaginary quadratic number field of class number $1$, let $\OO$ be its maximal
order, and let $\mathcal{F}_K$ be as defined above.

\begin{enumerate}
\item $\mathcal{F}_K$ is a fundamental domain for the action of
$\PGL_2(\OO)$ on $\mathcal{H}_3$. No two points in $\mathcal{F}_K$ are $\PGL_2(\mathcal{O})$-equivalent.
\item There exists a constant
$t_K>0$ such that $t\geq t_K$ for every $z+tj\in\mathcal{F}_K$. The value of
$t_K^2$ is given in the following tables : \vskip.5cm {\centering
\begin{tabular}{|c|c|c|c|c|c|} \hline
D & $1$ & $2$ & $3$ & $7$ & $11$   \\ \hline
\hline
$t_K^2$ & $1/2$ & $1/4$ & $2/3$ & $3/7$ & $2/11$ \\ \hline
\end{tabular}\par}
\vskip.5cm
{\centering \begin{tabular}{|c|c|c|c|c|} \hline
D & $19$ & $43$ & $67$ & $163$   \\ \hline
\hline
$t_K^2$ & $2/19$ & $2/43$ & $2/67$ &  $2/163$   \\ \hline
\end{tabular}\par}
\end{enumerate}
\end{theo}
\preuve
\begin{enumerate}
\item Since $\PGL_2(\OO) /\PSL_2(\OO)\simeq \OO^{\times}/(\OO^{\times})^2$ our hypotheses imply that its cardinality is $2$. Using the well-known fundamental domains for the 
  $\PSL_2(\mathcal{O})$ action on $\mathcal{H}_3$ (see for example \cite{EGM}) the result follows.

The action of $\PGL_2(\mathcal{O})$ on the boundary of $\mathcal{H}_3$ is generated by the following matrices:
\begin{enumerate}
\item Either  $\left(\begin{array}{cc}0 & -1 \\ 1 & 0
      \end{array}\right)$ acting on points $z+jt$ such that
  $\abs{z}^2+t^2=1$.
\item Or translations of the form  $\left(\begin{array}{cc}1 & \alpha \\ 0 & 1 \end{array}\right)$,
  for an appropriate $\alpha\in\mathcal{O}$.
\end{enumerate}
A tedious computation yields the result.
\item See \cite{Cre0} and \cite{Wit} for details.\end{enumerate} \eproof

\misctitle{Remark}
Thanks to Definition \ref{fund_dom} and Theorem \ref{boundt}, we have
explicit bounds for $z$ and $t$-components of elements in a fundamental
domain of $\mathcal{H}_3$ modulo $\GL_2(\OO)$, when $\OO$ is principal.
Unfortunately, when $h_K\neq 1$, we do not have a lower bound for $t$ (there
are points in the boundary of the fundamental domain such that $t=0$), and this will prevent
us from bounding the coefficients of reduced forms. This is the reason why we
will restrict our work to the class number $1$ case.

\subsection{Reduction of binary Hermitian forms}
Before tackling cubic forms, we recall the classical reduction theory of
binary Hermitian forms modulo $\GL_2(\OO)$, where $\OO$ is the maximal
order of an imaginary quadratic field.
\begin{deff}
Let $(P,Q,R)$ denote the binary Hermitian form 
$$H(x,y)=P\vert
x\vert^2+Q\overline{x}y+\overline{Q}x\overline{y}+ R\vert y\vert^2,
\quad P,Q,R\in \mathbb{C},$$ of
discriminant $\disc(H) :=-\Delta =  \vert Q\vert ^2 - PR$, and let $\mathscr{P}$ be the
set of \emph{positive definite} binary quadratic Hermitian forms over
$\mathbb{C}$; in other words,
$$ \mathscr{P} = 
 \left\{(P,Q,R) \colon P,R\in\mathbb{R}^+, Q\in\mathbb{C}, \disc(P,Q,R) <0\right\}.$$
The group $\PGL_2(\mathcal{O})$ acts on $\mathscr{P}$ via
$$M\cdot H(x,y)=H(\alpha x+\beta y,\gamma x+\delta y),
\quad
\text{where $M=\smallmat{\alpha & \beta \\ \gamma & \delta}\in \PGL_2(\mathcal{O})$}.
$$
\end{deff}

\misctitle{Remark}
It is customary to identify the Hermitian form
$$
\begin{pmatrix}
\overline{x} & \overline{y}
\end{pmatrix}
\begin{pmatrix} P & Q\\
\overline{Q} & R
\end{pmatrix}
\begin{pmatrix} x \\ y \end{pmatrix}
$$
with the Hermitian matrix $H = \smallmat{P & Q \\ \overline{Q} & R};$
the $\PGL_2(\mathcal{O})$ action is then
$$M\cdot H =M^*\times H \times M,
\quad\text{where $M^*=(\overline{M})^t$}.$$

\begin{lem}
Let $\widetilde{\Phi}: \mathscr{P}/\mathbb{R}_+^*\rightarrow \mathcal{H}_3$ be defined by:
\begin{equation}\label{Phi}\Phi\big( (P,Q,R) \big)
=-\frac{Q}{P}+\frac{\sqrt{\Delta}}{P}j.\end{equation}
$\widetilde\Phi$ is a bijection wich commutes with the action of $\PGL_2(\mathcal{O})$.
\end{lem}

This defines natural representatives for orbits of Hermitian
forms modulo $\PGL_2(\mathcal{O})$. Namely

\begin{deff}\label{hred} Let  $H\in\mathscr{P}$ a binary hermitian form. $H$ is called
\emph{reduced} if and only if $\widetilde\Phi(H)\in \mathcal{F}_K$.
\end{deff}

\begin{lem}\label{lemmaPQR}
Let
$(P,Q,R)=P\abs{x}^2+Qx\overline{y}+\overline{Q}\overline{x}y+R\abs{y}^2$
be a reduced Hermitian form in $\mathscr{P}$, with discriminant
$-\Delta = \abs{Q}^2 - PR$.
We have 
\begin{equation}\label{1ineq}P\leq \frac{\sqrt{\Delta}}{t_K}.\end{equation}

\begin{equation}\label{2ineq}\vert Q\vert ^2\leq c_K P^2,\end{equation}
and
\begin{equation}\label{3ineq}PR\leq \left(1+\frac{c_K}{t_K^2}\right)\Delta,\end{equation}
where $c_K$ is a constant depending only on the number field $K$, defined as follows 
$$c_k=\left\{\begin{array}{ll}1/2 & \textrm{if }K=\mathbb{Q}(i)\\ 7/12 & \textrm{if }K=\mathbb{Q}(\sqrt{-3})\\ \left(\frac{1+\abs{d_K}}{4}\right) & \textrm{otherwise} \end{array}\right. .$$\end{lem}
\preuve

\noindent For \eqref{1ineq} just recall that $t=\sqrt{\Delta}/P$ by the
definition of $\widetilde\Phi$ in \eqref{Phi} and $t\geq t_K$.

Thanks to the bounds on $\re{(z)}$ and $\im{(z)}$ given in the description
of the fundamental domain $\mathcal{F}_K$ (in Definition \ref{fund_dom}) we get
\begin{itemize}
\item $0\leq\vert \re(Q)\vert\leq P/2$, $0\leq \im(-Q)\leq 1/2$, and
  so $\vert Q\vert^2\leq P^2/2$ when $K=\mathbb{Q}(i)$;
\item $0\leq \re(-Q)\leq P/2$, $-\sqrt{3}/6P\leq \im(-Q)\leq \sqrt{3}/3P
  $ and then $\vert Q\vert^2\leq 7/12 P^2$, when $K=\mathbb{Q}(\sqrt{-3})$;
\item 
$0\leq  \re(-Q)\leq P/2$,
$0\leq  \im(-Q)\leq \frac{\sqrt{\vert d_K\vert}}{2}P$ and then
$  \vert Q\vert^2\leq \left(\frac{1+\vert
    d_K\vert}{4}\right)P^2$.
\end{itemize}
In all cases we have
$$\vert Q\vert^2\leq c_K P^2\leq c_K\frac{\Delta}{t_K^2}.$$
Recalling that $PR-\vert Q\vert^2 =  \Delta$, we obtain
$$PR\leq \left(1+\frac{c_K}{t_K^2}\right)\Delta.$$
\eproof

\subsection{Julia's covariant}
From now on, let $K$ be an imaginary quadratic field, let $\OO$ be its ring
of integers, and let $ V_{\mathcal{O}}$ be the set of binary cubic forms in
$\mathcal{O}[x,y]$. We want to define a canonical representative (or
\emph{reduced form}) in each orbit $\GL_2(\mathcal{O}) \cdot F$, $F\in
V_{\mathcal{O}}$.

\begin{deff}
We consider binary cubic forms in $V_{\mathcal{O}}$,
$$F(x,y)=ax^3+bx^2y+cxy^2+dy^3,\quad a,b,c,d\in\mathcal{O}$$
modulo the action of $\GL_2(\mathcal{O})$ given by
 $$M\cdot F=(\det(M))^{-1}F(Ax+By,Cx+Dy),\quad \textrm{ for each }M=\begin{pmatrix}A & B\\C &
    D\end{pmatrix}\in\GL_2(\mathcal{O}).$$
\end{deff}
\misctitle{Remark} As we saw in Corollary \ref{hk1}, this is the restriction
of the action used in Taniguchi's Theorem, when $h_K=1$. \vskip.5cm

\noindent Julia~\cite{Jul} gives us a covariant for this action:
\begin{deff}
Let $F\in V_{\mathcal{O}}$ be irreducible over $K$, factoring over $\mathbb{C}$
as $F(x,y) = a(x-\alpha_1 y)(x-\alpha_2 y)(x-\alpha_3 y)$, with $a\neq
0$. We associate to $F$ the positive definite binary Hermitian form
$$H_F(x,y)=t_1^2\vert x-\alpha_1 y\vert^2+t_2^2\vert x-\alpha_2
y\vert^2+t_3^2\vert x-\alpha_3 y\vert^2,$$
where 
$$t_i^2=\vert a\vert^2\vert\alpha_j-\alpha_k\vert^2,\quad i,j,k
\textrm{ pairwise distinct.}$$
\end{deff}
\noindent The following three lemmas follow from a direct computation:
\begin{lem}
We have
$$H_F(x,y)=P\vert x\vert^2+Q
\overline{x}y+\overline{Q}x\overline{y}+R\vert y\vert^2,$$
where 
$$\left\{\begin{array}{l}P=t_1^2+t_2^2+t_3^2\in\mathbb{R}^{+},\\
  Q=-(\alpha_1 t_1^2+\alpha_2 t_2^2+\alpha_3 t_3^2)\in\mathbb{C},\\
  R=\vert
  \alpha_1\vert^2 t_1^2+\vert \alpha_2\vert^2 t_2^2+\vert
  \alpha_3\vert^2 t_3^2\in\mathbb{R}^{+}.\end{array}\right.$$
\end{lem}
\begin{lem}
We have
\begin{equation}
(t_1t_2t_3)^2=\vert a\vert^2\vert \disc(F)\vert
\end{equation}
\end{lem}
\begin{lem}
Let $\Delta=-\disc(H_F)=PR-\vert Q\vert^2$ and $D=\disc(F)$. Then
\begin{equation}\Delta=3\vert D\vert .\end{equation}
\end{lem}

\begin{prop}
The application which sends $F$ to $H_F$ is covariant, i.e.
$$H_{M\cdot F}=M\cdot H_F,$$
for all $M\in \GL_2(\mathcal{O})$.
\end{prop}

\noindent Thanks to this property we can translate our problem of defining a
unique reduced $F$ to the problem of finding a unique reduced
covariant $H_F$ plus some extra conditions as we will see in Section
\ref{aut}.  

\begin{deff}[Julia reduction]
Let $F=(a,b,c,d)\in V_{\mathcal{O}}$ be a binary cubic form
with coefficients in $\mathcal{O}$.
We say that $F$ is Julia-reduced (modulo $\GL_2(\mathcal{O})$) if its covariant $H_F$ is reduced, in the sense of Definition \ref{hred}.
\end{deff}
\section{Reduction of binary cubic forms}\label{sec_bounds}

\subsection{Bounds for binary cubic forms}\label{t_red}\label{boundsabcd}
Let $F$ be a binary cubic form and let $H_F$ be its covariant hermitian form.
Starting from bounds on $H_F$ coefficients it is possible to directly bound $F$ coefficients and then to loop over all reduced binary cubic forms in time $\widetilde{O}(X)$, but the coefficients involved in the complexity of this algorithm are quite  big (see \cite{Mor} for details), so we chose another method, suggested by John Cremona, which can be found in \cite{Cre,Wom,CreSlides}.

\begin{deff}
For any  $k\in\mathcal{O}$, we note  $\tau_k = \left(\begin{smallmatrix}1 & k \\ 0 &
  1\end{smallmatrix}\right)$.
\end{deff}

\begin{deff}
For any $a_0\in\mathcal{O}$, we fix once and for all a system of representatives $\mathcal{P}_{a_0}$ for $\mathcal{O}/3a_0\mathcal{O}$. This is a finite set with $9\abs{a_0}^2$ elements.

\end{deff}

\smallskip

\begin{deff}
Let $F_K$ be as in definition \ref{fund_dom}.

We define $\mathcal{P}_K$ to be a fundamental region for $\mathbb{C}/\mathcal{O}$ such that $F_K\subset\mathcal{P}_K$.

\end{deff}

\begin{prop}\label{boundsabcd}
Let $F=(a,b,c,d)$ be a binary cubic form. There exists a unique $k \in \mathcal{O}$ such that $\tau_k$ sends $F$ to an equivalent binary cubic form $F_0=(a_0,b_0,c_0,d_0)$ such that $b_0\in \mathcal{P}_{a_0}$. We will call this $F_0$ $\tau$-reduced.

Moreover, if $F$ is Julia-reduced, then we have also the following properties:

\begin{equation}\begin{array}{l}\abs{a_0}\leq 3^{-3/4}t_K^{-3/2}D^{1/4}\\
\vert c_0\vert\leq \frac{\vert b_0\vert^2 + c_H
  D^{1/2}}{3\vert a_0\vert} \\
\textrm{either} \abs{d_0-x_1}\leq \frac{X^{1/4}}{\sqrt{\abs{A}}} \textrm{ or } \abs{d_0-x_2}\leq \frac{X^{1/4}}{\sqrt{\abs{A}}},\end{array} 
\end{equation}
where $c_H = 3^{1/2} 2^{-1/3} t_K^{-1}$,  $A=-27a_0^2, B=18a_0b_0c_0-4b_0^3$,
$C=b_0^2c_0^2-4a_0c_0^3,$ and we call $x_1$ and $x_2$ the roots of the quadratic polynomial $Ax^2+Bx+C$.
\end{prop}

\preuve
\noindent As regards the first assertion, just remark that $\tau_k$ sends $(a,b,c,d)$ to $(a_0,b_0,c_0,d_0)=(a,b+3ak,3ak^2+2bk+c,ak^3+bk^2+ck+d)$.

Now, assume that $F$ is Julia reduced. 

Let us consider the seminvariants associated to $F_0$:
$$P_H = b_0^2-3a_0c_0\quad \textrm{and}\quad U_H=2b_0^3+27a_0^2d_0-9a_0b_0c_0$$
(recall that $P_H$ is the first coefficient of the Hessian of $F_0$, but it is not in general equal to $P_0$, the first coefficient of the covariant associated to $F_0$).

$\tau_k$ leaves $P_H$ and $U_H$ unchanged and, as shown
in Womack's thesis \cite{Wom} we have 
\begin{equation}\abs{a_0}\leq 3^{-3/4}t_K^{-3/2}X^{1/8}\label{bounda}\end{equation}
and
$$\abs{U_H}\leq 3^{3/4}t_K^{-3/2}X^{3/8}$$ 
so from the syzygy
$$4P_H^3=U_H^2+27 \disc(F_0) a_0^2$$
we obtain

\begin{equation}\label{ph}P_H\leq c_H X^{1/4},\end{equation}
where $c_H=3^{1/2} 2^{-1/3} t_K^{-1}$,
and we easily obtain the bound for $\abs{c_0}$.

Finally, since $\disc(F_0)=Ad_0^2+Bd_0+C$ and $\abs{\disc(F_0)}\leq \sqrt{X}$ we have
$$\abs{d_0-x_1}\abs{d_0-x_2}\leq \sqrt{X}/\abs{A},$$
and this inequality implies that $\abs{d_0-x_1}$ and $\abs{d_0-x_2}$ can't be both bigger than $\frac{X^{1/4}}{\sqrt{\abs{A}}}.$
\eproof

\begin{coro}\label{coro_time}
It is possible to list all the reduced binary cubic forms $(a,b,c,d)$ (modulo
$\GL_2(\mathcal{O})$), with $\mathcal{N}(\disc(F))\leq X$ in time $O(X^{1+\varepsilon} ),$
for all $\varepsilon>0$.
\end{coro}

\preuve

\noindent The number of $\tau$-reduced binary cubic forms $(a_0,b_0,c_0,d_0)$ which are equivalent to  Julia-reduced ones (i.e. satisfying all properties enumerated in Proposition \ref{boundsabcd}) is 
$$N \ll\sum_{\abs{a_0} \ll X^{1/8}}\,\sum_{b_0\in\mathcal{P}_{a_0}}\,\sum_{c_0\ll X^{1/4}/\abs{a_0}}\, \sum_{\substack{\abs{d_0-x_1}\ll X^{1/4}/\abs{a_0} \\ \textrm{or} \abs{d_0-x_2}\ll X^{1/4}/\abs{a_0}}}1$$
Thus
\begin{equation*}N \ll  \sum_{\abs{a_0}\ll
    X^{1/8}}\abs{a_0}^2\cdot \frac{X^{1/2}}{\abs{a_0}^2}\cdot \frac{X^{1/2}}{\abs{a_0}^2} = X\cdot\sum_{\abs{a_0}\ll X^{1/8}}\frac{1}{\abs{a_0}^2}\end{equation*}
And

$$\sum_{\abs{a_0}\ll
X^{1/8}}\frac{1}{\abs{a_0}^2}
 \ll \sum_{n=1}^{X^{1/4}}
  \frac{\#\{a_0\in\mathcal{O}\colon
  \abs{a_0}^2=n\}}{n}
.$$
Since $\#\{a_0\in\mathcal{O}\colon
  \abs{a_0}^2=n\}=O(n^\varepsilon) = O(X^{\varepsilon})$ for all $\varepsilon>0$, and
  $\sum_{n=1}^{X^{1/4}}\frac{1}{n}$ is $O(\log (D)),$ we can conclude.
\eproof

\begin{prop}\label{backjulia}
Let $F$ be a Julia-reduced binary cubic form, let $F_0$ be the corresponding $\tau$-reduced form, and let $H_{F_0}=(P_0,Q_0,R_0)$ be the binary Hermitian form associated to $F_0$. Then $F=\tau_{-k}\cdot F_0$, for a unique $k\in\mathcal{O}$.
\end{prop}
\preuve

The action of $\tau_k$ sends $H_{F} = (P,Q,R)$ to $H_{F_0}=(P_0,Q_0,R_0)$ such that $P_0=P$ and $Q_0=Q+kP$. Dividing by $P$, we obtain $z_0=z-k$, with $z\in F_K\subset \mathcal{P}_K$, but this uniquely determines $k$, so we can conclude.
\eproof
\begin{algo}[$\tau$-reduction]\label{taured}
Let $K$ be an imaginary quadratic number field of class number $1$.
This algorithm loops over all $\tau$-reduced binary cubic forms $F_0=(a_0,b_0,c_0,d_0)$ satisfying the conditions in Proposition \ref{boundsabcd} with
$\mathcal{N}(\disc(F_0)) \leq X$, and associates the equivalent binary cubic form $F=(a,b,c,d)=\tau_{-k}\cdot F_0$, such that $z\in F_K$ (as explained in Proposition \ref{backjulia}).  
\vskip.5cm
\sf
For each $a_0,b_0,c_0,d_0$ in $\mathcal{O}$ satisfying the following properties:
\begin{itemize}
  \item $\abs{a_0} \leq \left(\frac{1}{t_K\sqrt{3}}\right)^{3/2}X^{1/8}$,
  \item  $b_0$ belongs to $\mathcal{P}_{a_0}$
  \item $\vert c_0\vert \leq \frac{\vert b_0\vert^2 + c_HX^{1/4}}
                         {3\vert a_0\vert}, $
  \item Either $\abs{d_0-x_1}\leq X^{1/4}/\sqrt{\abs{A}}$ or $\abs{d_0-x_2}\leq X^{1/4}/\sqrt{\abs{A}} $
\end{itemize}
Do the following operations:
\begin{enumerate}
\item compute the first two coefficients $P_0,Q_0$ of the covariant $H_{F_0}$
of the cubic form $F_0=(a_0,b_0,c_0,d_0)$.
\item Compute $k$ such that $z_0+k\in \mathcal{P}_K$ ($z_0=-Q_0/P_0)$.
\item Compute $F=(a,b,c,d)=\tau_{-k}(a_0,b_0,c_0,d_0).$
\end{enumerate}
\end{algo}

\subsection{Automorphism matrices}\label{aut}

In this section we are going to study automorphism matrices for binary hermitian forms.

\begin{prop}\label{autbounds}
Let $F=(a,b,c,d)$, $F$ Julia-reduced.
Let $H=H_F$, and $\Delta = PR-\abs{Q}^2$.
 
Let $M=\left(\begin{array}{cc}A & B \\ C & D\end{array}\right)\in
  \GL_2(\mathcal{O})$ such that $M\cdot H=H$.
Then we have the following bounds on the coefficients of $M$:
\begin{equation}\abs{A}^2\leq \frac{PR}{\Delta},\quad\abs{C}\leq
  \frac{P}{\sqrt{\Delta}},\quad \abs{D}^2\leq
  \frac{PR}{\Delta}.\end{equation}
and
\begin{equation}\left\{\begin{array}{ll}\abs{B}\leq \frac{PR}{\Delta}+1 & \textrm{if
  }B\neq 0\\ \abs{B}\leq 2\sqrt{c_K}& {\textrm{if }B=0} \end{array}\right..\end{equation}
\end{prop}
\preuve

\noindent Let us write $H(x,y)=P\vert
x\vert^2+Q\overline{x}y+\overline{Q}x\overline{y}+R\vert y\vert^2$.
We have
\begin{equation}\label{PH}PH(x,y)=\vert xP+yQ\vert^2+\Delta\vert
  y\vert^2,\quad \textrm{and}\end{equation}
\begin{equation}RH(x,y)=\vert Ry+\overline{Q}x\vert^2+\Delta\vert x\vert^2.\end{equation}
Thanks to formula (\ref{PH}) we can give upper bounds for $\vert
A\vert,\vert B\vert$,  and $\vert D\vert$. Let us write more explicitly the
relation $M\cdot H = H$:
\begin{eqnarray}\label{automat}M\cdot H&=&\left(\begin{array}{lr}\overline{A}&\overline{C}
  \\ \overline{B}&\overline{D}\end{array}\right)\left(\begin{array}{lr}P & Q\\ \overline{Q}
  &R\end{array}\right)\left(\begin{array}{lr}A & B\\C &
    D\end{array}\right)
\\&=&\left(\begin{array}{lr}\vert A\vert^2
    P+\overline{A}CQ+A\overline{C}\overline{Q}+\vert C\vert^2R &
    \overline{A}BP+\overline{A}DQ+B\overline{C}\overline{Q}+\overline{C}DR\\
   A\overline{B}P+C\overline{B}Q+A\overline{D}\overline{Q}+C\overline{D}R
    & \vert B\vert^2P+\overline{B}DQ+B\overline{D}\overline{Q}+\vert
    D\vert^2 R\end{array}\right)\nonumber\\ &= &\left(\begin{array}{cc}H(A,C)
    & \dots \\ \dots & H(B,D)\end{array}\right).\nonumber\end{eqnarray}
By imposing this matrix to be equal to $M$ we have
\begin{eqnarray*}
\vert AP+CQ\vert^2+\Delta\vert C\vert^2=P^2 &\Rightarrow& \vert C\vert\leq
\frac{P}{\sqrt{\Delta}},\\
\vert BP+DQ\vert^2+\Delta\vert D\vert^2=PR &\Rightarrow & \vert D\vert^2 \leq
\frac{PR}{\Delta},\\
\vert CR+A\overline{Q}\vert^2+\Delta\vert A\vert^2=PR &\Rightarrow& \vert
A\vert^2 \leq \frac{PR}{\Delta}.\end{eqnarray*}
When $C=0$ the third equation becomes 
$$A\overline{B}P+A\overline{D}\overline{Q} = \overline{Q},$$ with $\abs{A}=\abs{D}=1$ so it
is easy to check that $\abs{A\overline{C}}P\leq 2 Q\leq 2\sqrt{c_K}P$
and we obtain the formula.
Finally, when $C\neq 0$, since $\abs{AD-BC}=1$ we get $$\abs{B}\leq
\frac{1+\abs{AD}}{\abs{C}}$$ and we easily conclude.
\eproof

\noindent The bounds of the previous Proposition are completely explicit when
$h_K=1$, since we know $t_K$ and $c_K$.

\begin{deff}
Let $M\in\PGL_2(\mathcal{O})$. We define 
$$S(M)=\{H\in \mathscr{P}/\mathbb{R}_+^* \mid M\cdot H = H \textrm{ and }H\textrm{ reduced}\}.$$
that is, the set of reduced binary hermitian forms which are stabilized by the action of $M$.
\end{deff}

\medskip

\noindent The
following algorithm lists the finite set of automorphism matrices. It needs to be run only once for each of our 9
imaginary quadratic fields of class number~1.
\begin{algo}\label{al:autmat}
Computes the set $\mathcal{M}$ of all matrices $M$ stabilizing some reduced binary hermitian form, and for each $M$ outputs also the corresponding set $S(M)$.
\sf
\medskip

Set $\mathcal{M}=\emptyset$.

For each triple $(A,C,D)$ satisfying the bounds of Proposition \ref{autbounds}, do the following operations:
\begin{enumerate}
\item For each $B\in \mathcal{O}$ such that $\vert AD-BC\vert = 1$ (if $C=0$,
take only the set $\abs{B}\leq 2\sqrt{c_K})$, let  $M=\left(\begin{array}{cc}A & B \\ C & D\end{array}\right)$ and do the following.

\item Consider the following $4\times 4$ matrix, with coefficients in
$\mathcal{O}$:
$$W(M)=\left(\begin{array}{cccc}(\vert A\vert^2-1) & \overline{A}C &
    A\overline{C} & \vert C\vert ^2\\
\overline{A}B & (\overline{A}D-1) & B\overline{C} & \overline{C}D\\
A\overline{B} & C\overline{B} & (A\overline{D}-1) & C\overline{D} \\
\vert B\vert^2 & \overline{B}D & B\overline{D} & (\vert
D\vert^2-1)\end{array}\right).$$
\item Compute the rank $r$ of $W(M)$ (over the field $K$).
\item {\bf If} $r=1$ or $r=4$, {\bf skip} to the following quadruple $(A,B,C,D)$.
\item {\bf If} $r=0$ {\bf output } $M=\left(\begin{array}{cc}A & B\\C
  & D\end{array}\right)$ and $S(M)=\{H\in\mathscr{P}/\mathbb{R}_+^*\mid H \textrm{ reduced}\}$. $\mathcal{M}=\mathcal{M}\cup\{M\}$.
\item {\bf If} $r=2$ or $r=3$, set  $M=\left(\begin{array}{cc}A & B\\C
  & D\end{array}\right)$, compute the set $S(M)=\{H= (P,Q,R)\in \mathscr{P}/\mathbb{R}^*_+\mid W\cdot (P,Q,\overline{Q},R)^t = 0 \textrm{ and } (P,Q,R) \textrm{ reduced}\}$. If $S(M)\neq
\emptyset$, {\bf output } $M=\left(\begin{array}{cc}A & B\\C
  & D\end{array}\right)$ and $S(M)$.  $\mathcal{M}=\mathcal{M}\cup\{M\}$.
\end{enumerate}
Output $\mathcal{M}$.
\end{algo}
\misctitle{Remarks}
\begin{itemize}

\item We could also loop only on $A,D$ and replace step (1) by :
\begin{enumerate}
\item Solve $\vert AD-BC\vert  = 1$ for $B,C\in\mathcal{O}$. This time $BC$
belongs to an explicit finite set, and we enumerate divisors.
\end{enumerate}
\item It is possible to write (once for all) explicit conditions to associate to any binary hermitian form $H$ its set of automorphism matrices, just looking at the sets $S(M)$ computed in algorithm \ref{al:autmat}
\item For an example of application of the above algorithm, Appendix \ref{app_aut} contains the list of all automorphism matrices for $K=\mathbb{Q}(i)$ and the corresponding conditions on binary hermitian forms.
\end{itemize}

\misctitle{Remark}
Running the algorithm on all the $9$ possible number fields we noticed a property holding for all $K\neq \mathbb{Q}(\sqrt{-3})$ :
\begin{itemize}
\item For each matrix $M$ found at step 6 (that is, they are not trivial automorphisms) $W(M)$ has rank $2$ and $S(M)$ is a subset of the boundary of the fundamental domain.
\end{itemize}
In the case $K=\mathbb{Q}(\sqrt{-3})$ we have explicit counter-examples.

\vskip.5cm
\noindent The proof of Algorithm \ref{al:autmat} is given by the following proposition.
\begin{prop}\label{autmat}
Let $M = \smallmat{A&B\\C&D}\in \PGL_2(\mathcal{O})$
belong to the stabilizer of $H_F$, where $H_F$ is the Hessian of some reduced cubic
form $F$. If $r$ is the rank of the matrix $W$ constructed in the above
algorithm, then
\begin{itemize}
\item $r=0$ if and only if $B = C = 0$ and $A = D$ are units. Then $M$ is
an automorphism for \emph{all} Hermitian quadratic forms in $\mathcal{F}$.
\item $r=1$ is impossible
\item $r=2$ or $r=3$  then $M$ is an automorphism for some
  linear subspace of  $\mathscr{P}$, defined by explicit equations in
the variables $P,Q,\overline{Q},R$.
\item $r=4$ is impossible.
\end{itemize}
\end{prop}
\preuve

\noindent The condition $\smallmat{A&B\\C&D}\in \Aut(H)$ translates to the linear system $W(M)\cdot X = 0$, with 
$X = (P,Q,\overline{Q}, R)^t$.
\begin{itemize}
\item If $r=4$, the only solution of the system is $(0,0,0,0)$, but this is
not allowed since $P,R >0.$

\item Assume that $r\leq 1$ : the matrix $\smallmat{A & B \\ C & D}$ has rank
2 so the two $2$ by $2$ matrices on the lower-left and upper-right corners of
$W(M)$ have rank $2$ unless $B=C=0$. In this case $W(M)$ is diagonal
$$\left(\begin{array}{cccc}\vert A \vert^2 -1 & & &
  \\ & \overline{A}D-1 & &\\ & & A\overline{D}-1  & \\ & & & \vert
  D\vert^2-1\end{array}\right).$$
Since $B=C=0$, and $AD - BC$ is a unit, we must have $\vert A \vert
=\vert D\vert =1$, so this matrix has either rank $2$ or $0$ (when
$\overline{A}D=A\overline{D}=1$).
\end{itemize}
\eproof

\section{The algorithm}\label{sec_algo}

\begin{algo}\label{main_algo}
Given a bound $X=D^2$, output the list of reduced binary cubic forms
modulo $\GL_2(\mathcal{O})$, such that  $\mathcal{N}(\disc(F))\leq X$.
\vskip.2cm
\sf

\noindent Use sub-Algorithm \ref{taured} to loop over quadruples $F=(a,b,c,d)\in \mathcal{O}^4$ satisfying all the
properties in  Section \ref{boundsabcd}. Do the following operations:
\begin{enumerate}

\item Approximate the complex roots of $F$,
  $(\alpha_1,\alpha_2,\alpha_3)$ to a sufficient accuracy. Then approximate $H=H_F=(P,Q,R)$ the associated Hermitian form.
\item Check if $H$ is in the fundamental domain modulo
  $\PGL_2(\mathcal{O})$ (i.e. it is reduced), (see Definition
  \ref{fund_dom}). In particular, if $H_F$ is ``near'' to
  the boundary of the fundamental domain use Algorithm \ref{borders}
  (see below) to check exactly the boundary condition.
  If not {\bf skip} to the following $F$.

\item Check whether $F$ is irreducible in $K[x,y]$. If not
{\bf skip} to the following $F$.
\item Apply Dedekind criterion to check whether $F$ describes a maximal
  ring. If not {\bf skip} to the following $F$.
\item Apply sub-algorithm \ref{al:autmat} to compute $\mathcal{M}$, the set of all automorphism matrices for $H$.
\item Compute the set $\{M\cdot F\mid M\in\mathcal{M}\}$ and check if $F$ is the minimal element of this set (for some order, for instance the lexicographic one). If not {\bf skip} to the following $F$.
\item {\bf print} $F$. 

\end{enumerate}
\end{algo}

\misctitle{Remarks}
\begin{itemize}
\item For the precision needed in step $(1)$ refer to Appendix C of \cite{Mor}.
\item In step (5), we compute a list of automorphs for $F$
  to decide whether $F$ is minimal among the reduced forms in its orbit
  with respect to the lexicographic order
  (in this case $F$ should be kept, otherwise not). Another way to deal with
  this problem would be to store all those $F$ and then check
  $\GL_2(\mathcal{O})$-equivalence once we have all the forms with a
  fixed discriminant $D$. The problem is that our algorithm does not output
  forms ordered by discriminant, so we could apply this test only at the end, 
  and this would increase dramatically the space complexity. (Remember that we
  output the ``good'' binary cubic forms as we find them, so we do not keep in
  memory the list of representatives of cubic extensions).   
\end{itemize}

\section{Implementation problems}\label{sec_details}
\subsection{Checking rigorously the boundary conditions}\label{se:rigorous}
As the computation of $P,Q,R$ involves floating point approximations of the
complex roots of a polynomial in $\mathcal{O}[X]$, it will not give, of
course, exact results. Those floating point computations will in general be
sufficient to test whether the Hermitian form is strictly inside or outside
the fundamental domain. But if it is very near the boundary (or worse 
\emph{on} the boundary), this approach fails.

\noindent To get rid of this problem we use the following formulas:
\begin{eqnarray}
P&=&-\frac{\vert b\vert^2}{\vert a\vert^2}+3(\vert
\alpha_1\vert^2+\vert\alpha_2\vert^2+\vert\alpha_3\vert^2)\\
Q&=&\frac{\overline{b}c}{\vert a\vert^2}+3(\overline{\alpha_1}\alpha_2\alpha_3+\alpha_1\overline{\alpha_2}\alpha_3+\alpha_1\alpha_2\overline{\alpha_3})\\
R&=&-\frac{\vert c\vert^2}{\vert a\vert^2}+3(\vert\alpha_1\vert^2\vert\alpha_2\vert^2+\vert\alpha_1\vert^2\vert\alpha_3\vert^2+\vert\alpha_2\vert^2\vert\alpha_3\vert^2)
\end{eqnarray}
Now we consider
$\alpha_1,\alpha_2,\alpha_3,\overline{\alpha_1},\overline{\alpha_2},\overline{\alpha_3}$ as
algebraic numbers, and we let $S$ be the set of the six 
permutations fixing the $\alpha_i$, and acting as $S_3$
on the $\overline{\alpha_i}$.
The polynomial
$$g_P = \prod_{\sigma\in
  S}(X-\sigma(\alpha_1\overline{\alpha_1}+\alpha_2\overline{\alpha_2}+\alpha_3\overline{\alpha_3}))$$
vanishes at
  $\vert\alpha_1\vert^2+\vert\alpha_2\vert^2+\vert\alpha_3\vert^2,$
and its coefficients are symmetric in $(\alpha_1,\alpha_2,\alpha_3)$
and $(\overline{\alpha_1}, \overline{\alpha_2}, \overline{\alpha_3})$
independently. They can thus be expressed in terms of
$(b/a,c/a,d/a)$ and $(\overline{b/a},\overline{c/a},\overline{d/a)}$.
The polynomial $f_P(X) = g_P\left(\frac{X}{3} -
\frac{\abs{b}^2}{3\abs{a}^2}\right)$ vanishes at $P$ and belongs to $K[X]$.

In the same way we can compute polynomials in $K[X]$ vanishing at $Q$, $R$,
$\re(Q)$ or $\im(Q)$. Such polynomials are easily computed using a computer
algebra system like Maple (and it is sufficient to compute them once for
all).

We want to verify rigorously boundary conditions, for instance $P = R$: if
$f_P$ and $f_R$ have no common factor in $K[X]$, then $P\neq R$. But this is
not enough: we also want to check whether $P < R$ or $P > R$, i.e.~if the
point we are testing is ``inside'' or ``outside'' the fundamental domain.

The following theorem of Mahler~\cite{Mah} provides the accuracy we need 
for our floating point computations:
\begin{theo}[Mahler]\label{mahler}
Let  $f=a_0x^m+a_1x^{m-1}+\dots+a_m=a_0(x-\alpha_1)\cdots(x-\alpha_m)$ be a
separable polynomial of degree $m\geq 2$, and let
$$\Delta(f)=\min_{1\leq i<j\leq m}\vert \alpha_i-\alpha_j\vert$$
be the minimal distance between two distinct roots of $f$. Then
$$\Delta(f)>\sqrt{3}m^{-(m+2)/2}\vert \disc (f)\vert^{1/2}M(f)^{-(m-1)},$$
where $\disc(f)$ is the discriminant of $f$, and $M(f)=\vert a_0\vert
\prod_{h=1}^m \max(1,\vert \alpha_h\vert)$.
\end{theo}

\noindent This translates to the following algorithm:

\begin{algo}[Checking an algebraic identity]\label{borders}
Let $\alpha$ and $\beta\in \mathbb{R}$ be two algebraic numbers, and let $A$
and $B\in K[X]\setminus {0}$ that vanish at $\alpha$, and $\beta$
respectively. Assume we can compute floating point approximations
$\hat{\alpha}$ and $\hat{\beta}$ such that $\abs{\alpha - \hat{\alpha}} <
\varepsilon$, $\abs{\beta - \hat{\beta}} < \varepsilon$, for any fixed
$\varepsilon > 0$.

We want to decide whether $\alpha < \beta$, $\alpha > \beta$ or $\alpha =
\beta$.

\sf
\begin{enumerate}
\item Let $C = AB$ and $f = C/\text{gcd}(C,C')$.

\item If the degree of $f$ is $1$, then {\bf answer} $\alpha = \beta$.

\item Compute a good approximation $\hat{\Delta}$ of 
$$\Delta(f)=\sqrt{3} m^{-(m+2)/2}\vert \disc(f)\vert^{1/2}M(f)^{-(m-1)},$$
where $\disc(f)$ and $M(f)$ are
defined in Theorem \ref{mahler} such that $\hat{\Delta} \leq \Delta(f)$.
\item Compute $\alpha$ and $\beta$ at precision $\varepsilon=\hat{\Delta}/4$,
i.e. $\hat{\alpha}$ and $\hat{\beta}$ such that
$$\abs{\alpha - \hat{\alpha}} < \varepsilon,\quad
  \abs{\beta - \hat{\beta}} < \varepsilon.$$

\item\label{st:5} If $\vert \hat{\alpha}-\hat{\beta}\vert < 2\varepsilon$,
 {\bf answer} $\alpha = \beta$.
\item If $\hat{\alpha} < \hat{\beta}$, {\bf answer} $\alpha < \beta$.
\item If $\hat{\alpha} > \hat{\beta}$, {\bf answer} $\alpha > \beta$.
\end{enumerate}
\end{algo}
\preuve

\noindent The polynomial $f$ is non constant and has $\alpha$ and $\beta$ among its
roots. If its degree is $1$, then $\alpha = \beta$. Otherwise,
assume first that $\vert \hat{\alpha}-\hat{\beta}\vert < 2\varepsilon$.
Then 
$$\abs{\alpha - \beta} \leq
\abs{\alpha - \hat{\alpha}}
+ \abs{\beta - \hat{\beta}}
+ \abs{\hat{\alpha} - \hat{\beta}}
< 4\varepsilon \leq \Delta(f).
$$
Hence $\alpha = \beta$ by Mahler's theorem in this case, proving \eqref{st:5}.

We now assume that $\vert
\hat{\alpha}-\hat{\beta}\vert \geq 2\varepsilon$; since
$$\alpha - \beta = \hat{\alpha} - \hat{\beta} 
+ (\alpha - \hat{\alpha})
- (\beta - \hat{\beta})$$
and $$\abs{(\alpha - \hat{\alpha}) - (\beta - \hat{\beta})} < 2\varepsilon,$$
$\alpha - \beta$ and $\hat{\alpha} - \hat{\beta}$ have the same sign.
\eproof

\begin{prop}
The smallest $\varepsilon$ that we can obtain in step (4) of the above
algorithm (i.e. the maximal precision needed) is $\gg X^{-\beta},$ for
some positive constant $\beta$. 
\end{prop}

\misctitle{Remark}
That means that for our computation we
will need at most $\Omega(\log X)$ significant digits.
\vskip.5cm
\preuve

\noindent Algorithm \ref{main_algo} loops over reduced integral cubic forms $F = (a,b,c,d)\in V_{\mathcal{O}}$ with discriminant $\disc(F)$ satisfying $\mathcal{N}(\disc(F))
\leq X$. In particular, Proposition~\ref{boundsabcd} implies that $\abs{a}\ll
X^{1/8}$.

\noindent For each such form, we may compute various separable polynomials $f$ with
coefficients in $a^{-u}\mathcal{O}_K$, for some bounded integer $u$. Then
$\disc(f)$ is non zero, in $a^{-4u}\mathcal{O}_K$. Its norm is a non-zero
rational integer divided by $\abs{a}^{-8u}$, hence $\gg X^{-u}$. Thus
$\disc(f)\gg X^{-u/2}$.

Landau's theorem (see \cite[Proof of Theorem~13.1]{Bel3} for example) tells us
that
$$M(f)\leq \Vert f\Vert_2$$
and the coefficients of $f$ are monomials in
$e_1,e_2,e_3,f_1,f_2,f_3$ (see Appendix D of \cite{Mor}). Each one of these is bounded by $c\cdot
X^{\alpha}$, for an appropriate constant $c$ and exponent $\alpha$.

We have 
$$\Delta(f)\gg M(f)^{-(m-1)}.$$
So we obtain $$\Vert f\Vert_2 \ll X^{\beta},$$
but then we can conclude that $\Delta(f)\gg X^{-\beta}$.\eproof

\subsection{An idea to count only half of the extensions}

Let $K$ be an imaginary quadratic number field, with class number $h_K=1$
and discriminant $d_K\neq -3,-4$.
It is easy to remark that if $H=(P,Q,R)$ is in the fundamental domain,
then also $H'=(P,-\overline{Q},R)$ is. And, in general, these two
Hermitian forms are not equivalent modulo $\PGL_2(\mathcal{O})$.

In particular, if $F=(a,b,c,d)$ has $H_F=H$, then
$F'=(\overline{a},-\overline{b},\overline{c},-\overline{d})$ gives
$H_{F'}=H'$.

So we can loop only on half of the $c$ satisfying the given bounds, then construct both the forms $F=(a,b,c,d)$ and
$F'=(\overline{a},-\overline{b},\overline{c},-\overline{d})$ and
check if they are equivalent (comparing $F'$
with the list of automorphic functions to $F$). If not we verify also the list
of automorphic functions to  $F'$ to see if one of them will be found
in our
loops, and if both answers are no, we add this second form $F'$ to
our output list. 

\section{Results}\label{se:results}
In this section we present results obtained for the case $K=\mathbb{Q}(i)$ via an implementation of our algorithm in Pari/GP \cite{PARI}, running on an Intel Xeon 5160 dual core, 3.0 GHz.

Let $X$ be the bound on $\mathcal{N}(\mathfrak{d}(L/K))$ and $N(X)$ the number of isomorphism classes of cubic extensions of $\mathbb{Q}(i)$ up to that
bound.

In the following table, we will compare the time needed to list all $N(X)$ representatives of cubic extensions of $\mathbb{Q}(i)$ with two algorithms : ray class fields one and ours. $t$ denotes the running time of our algorithm, $t'$ the running time of  ray class field one (see Section 9.2.3 of \cite{Coh2}).

\vskip1cm
 {\centering\begin{tabular}{|c|c|c|c|}

\hline
$X$ & $N(X)$ & $t$ & $t'$ \\ \hline
\hline
$10^4$ & $276$ & 5 s & 16 s    \\ \hline
$4\cdot 10^4$ & $1339$ & 19 s & 1mn 18 s\\ \hline
$9\cdot 10^4$ & $3305$ & 56 s & 3mn 45 s \\ \hline
$10^6$ & $42692$ & 24 mn 1 s  & 2h 52mn 9 s \\ \hline
$4\cdot 10^6$ & $181944$ & 2 h 49 mn  & 34h 24 mn 8 s\\ \hline
$9\cdot 10^6$ & $421559$ & 9 h 37 mn&  $>$ 134 h\\ \hline
$10^8$ & $4990974$ & 359 h 25 mn  & $>$ 2720 h \\\hline
\end{tabular}\par}

\misctitle{Remarks}

\begin{itemize}
\item This computations allowed us to check the correctness of our results. In fact, we compared $N(X)$ up to $X=9\cdot 10^6$ with the results of the ray class field algorithm  and all results matched.
\item The last line of the table would have involved too long computations with ray class field algorithm, so we skipped it, and we give only a prevision on the running time needed.
\end{itemize}
\subsection{Roberts' conjecture and asymptotic predictions for $N(X)$}

Frank Thorne compared our numerical results with heuristic asymptotic
developments derived from the Datskovsky-Wright method \cite{DW}, in the
spirit of Roberts's conjecture (see \cite{TT}). Starting from Roberts' conjecture, Taniguchi and Thorne worked out the formula in the particular case when $k$ is an imaginary quadratic number field : 
$$N(X)=\frac{1}{12}\frac{\Res_{s=1}\zeta_K(s)}{\zeta_K(3)}X+\frac{1}{40}d_K^{-1/2}\Res_{s=1}\zeta_K(s)\frac{\sqrt{3}\Gamma(1/3)^6}{\pi^2}\frac{\zeta_K(1/3)}{\zeta_K(2)\zeta_K(5/3)}X^{5/6}.$$
The following table compares our values for $N(X)$ with Thorne's asymptotic data. The results are strikingly similar.
\vskip.5cm
 {\centering\begin{tabular}{|c|c|c|}

\hline
$X$ & $N(X)$ (Morra) & $N(X)$ (Taniguchi-Thorne) \\ \hline
\hline
$10^4$ & $276$ & $270.2$   \\ \hline
$10^6$ & $42692$ & $42655.6$ \\ \hline
$9\cdot 10^6$ & $421559$ & $421260$ \\ \hline
$10^8$ & $4990974$ &  $4990962$ \\ \hline
\end{tabular}\par}

\appendix

\section{Automorphism matrices for $\mathbb{Q}(i)$}\label{app_aut}
{\centering\begin{tabular}{|l|l|}

\hline
$M$ (modulo multiplication by $\varepsilon \in\{+1,-1,i,-i\}$) & conditions for $S(M)$  \\ \hline
\hline
$\left(\begin{array}{cc}1 & 0 \\ 0 & 1\end{array}\right)$ & for all $P,Q,R$  \\ \hline
$\left(\begin{array}{cc}0 & -1 \\ 1 & 0\end{array}\right)$ & $\textrm{if }P=R\textrm{ and }\re(Q)=0$   \\ \hline
$\left(\begin{array}{cc}0 & 1 \\ 1 & 0\end{array}\right)$ & $\textrm{if }P=R\textrm{ and }\im(Q)=0$   \\ \hline
$\left(\begin{array}{cc}0 & 1 \\ -i & 0\end{array}\right)$ & $\textrm{if }P=R\textrm{ and }\re(Q)=\im(Q)$   \\ \hline
$\left(\begin{array}{cc}0 & 1 \\ i & 0\end{array}\right)$ & $\textrm{if }P=R\textrm{ and }\re(Q)=-\im(Q)$   \\ \hline
$\left(\begin{array}{cc}0 & -1 \\ 1 & 1\end{array}\right),\left(\begin{array}{cc}1 & 1 \\ -1 & 0\end{array}\right)$ & $\textrm{if }P=R\textrm{ and }\re(Q)=P/2$   \\ \hline
$\left(\begin{array}{cc}0 & -1 \\ 1 & -1\end{array}\right),\left(\begin{array}{cc}-1 & 1 \\ -1 & 0\end{array}\right)$ & $\textrm{if }P=R\textrm{ and }\re(Q)=-P/2 $  \\ \hline
$\left(\begin{array}{cc}0 & 1 \\ 1 & i\end{array}\right),\left(\begin{array}{cc}-i & 1 \\ 1 & 0\end{array}\right)$ & $\textrm{if }P=R\textrm{ and }\im(Q)=P/2 $  \\ \hline
$\left(\begin{array}{cc}-1 & 0 \\ 1 & 1\end{array}\right)$ &$ \textrm{if }P=R, \re(Q)=P/2\textrm{ and }\im(Q)=0  $ \\ \hline
$\begin{array}{c}{\left(\begin{array}{cc}i & 0 \\ -i & 1\end{array}\right),\left(\begin{array}{cc}-1 & -1-i \\ 1 & 1\end{array}\right),  \left(\begin{array}{cc}-1 & -1-i \\ -i & 1\end{array}\right),} \\ \left(\begin{array}{cc}-i & 0 \\ 1 & 1\end{array}\right), \left(\begin{array}{cc}0 & i \\ 1 & -i+1\end{array}\right), \left(\begin{array}{cc}1 & 1 \\ -i & -i-1\end{array}\right),\\ \left(\begin{array}{cc}1 & -i+1 \\-i-1 & -1 \end{array}\right), \left(\begin{array}{cc}1 & 0 \\ -i-1 & -1\end{array}\right),\\ \left(\begin{array}{cc}1 & 1 \\ -i-1 & -1\end{array}\right),  \left(\begin{array}{cc}1 & -i \\ -i-1 & -1\end{array}\right),  \left(\begin{array}{cc}1 & -i \\ -1 & i-1\end{array}\right)\end{array}$ &$\begin{array}{l} \textrm{if }P=R, \re(Q)=P/2 \\ \\  \textrm{and }\im(Q)=P/2\end{array}$   \\ \hline
$\left(\begin{array}{cc}-1 & 0 \\ -1 & 1\end{array}\right)$ & $\textrm{if }P=R,\re(Q)=-P/2\textrm{ and }\im(Q)=0$   \\ \hline
$\begin{array}{c}\left(\begin{array}{cc}-i & 0 \\ -i & 1\end{array}\right),\left(\begin{array}{cc}-1 & 1-i \\ -1 & 1\end{array}\right), \left(\begin{array}{cc}-1 & 1-i \\ -i & 1\end{array}\right),\\ \left(\begin{array}{cc}i & 0 \\ -1 & 1\end{array}\right), \left(\begin{array}{cc}0 & -i \\ 1 & -i-1 \end{array}\right), \left(\begin{array}{cc}1 & -i \\ 1 & -i-1\end{array}\right),\\ \left(\begin{array}{cc}1 & -1-i \\ 1-i & -1\end{array}\right), \left(\begin{array}{cc}1 & 0 \\ 1-i & -1\end{array}\right), \\ \left(\begin{array}{cc}1 & -i \\ 1-i & -1\end{array}\right),\left(\begin{array}{cc}1 & -1 \\ 1-i & -1\end{array}\right),\left(\begin{array}{cc}1 & -1 \\ -i & i-1\end{array}\right)\end{array}$ &$\begin{array}{l} \textrm{if }P=R, \re(Q)=-P/2 \\ \\ \textrm{and }\im(Q)=P/2\end{array}  $ \\ \hline
$\left(\begin{array}{cc}-1 & 0 \\ -i & 1\end{array}\right)$ & $\textrm{if }P=R,\re(Q)=0\textrm{ and }\im(Q)=P/2 $  \\ \hline
$\left(\begin{array}{cc}1 & 0 \\ 0 & -1\end{array}\right),\left(\begin{array}{cc}1 & 0 \\ 0 & -i\end{array}\right),\left(\begin{array}{cc}1 & 0 \\ 0 & i\end{array}\right)$ & $\textrm{if }Q=0 $\\ \hline
$\left(\begin{array}{cc}1 & 1 \\ 0 & -1\end{array}\right)$ &$ \textrm{if }\re(Q)=P/2\textrm{ and }\im(Q)=0 $  \\ \hline
$\left(\begin{array}{cc}1 & i \\ 0 & -i\end{array}\right),\left(\begin{array}{cc}1 & 1+i \\ 0 & -1\end{array}\right),\left(\begin{array}{cc}1 & 1 \\ 0 & i\end{array}\right)$ &$ \textrm{if }\re(Q)=P/2 \textrm{ and }\im(Q)=P/2$ \\ \hline
$\left(\begin{array}{cc}1 & -1 \\ 0 & -1\end{array}\right)$ &$ \textrm{if }\re(Q)=-P/2\textrm{ and }\im(Q)=0$   \\ \hline
$\left(\begin{array}{cc}1 & i \\ 0 & i\end{array}\right),\left(\begin{array}{cc}1 & -1 \\ 0 & -i\end{array}\right),\left(\begin{array}{cc}1 & -1+i \\ 0 & -1\end{array}\right)$ &$ \textrm{if }\re(Q)=-P/2 \textrm{ and }\im(Q)=P/2$ \\ \hline
$\left(\begin{array}{cc}1 & i \\ 0 & -1\end{array}\right)$ &$ \textrm{if }\re(Q)=0\textrm{ and }\im(Q)=P/2  $ \\ \hline
\end{tabular}\par}

\end{document}